\documentclass[12pt]{amsart}
\usepackage{amssymb,amsmath}
\usepackage{epsfig}
\usepackage{graphicx}
\usepackage{colordvi}
\newtheorem{theorem}[equation]{Theorem}

\newtheorem{proposition}[equation]{Proposition}

\theoremstyle{definition}

\newcounter{FNC}[page]
\def\fauxfootnote#1{{\addtocounter{FNC}{2}$^\fnsymbol{FNC}$%
     \let\thefootnote\relax\footnotetext{$^\fnsymbol{FNC}$#1}}}

\newcommand{\C}{{\mathbb C}}
\newcommand{\R}{{\mathbb R}}
\newcommand{\frakg}{{\mathfrak g}}

\begin{document}
\title[General isotropic flags are general]{General isotropic flags are general\\ (for
  Grassmannian Schubert calculus)}
\author{Frank Sottile}
\address{Department of Mathematics\\
         Texas A\&M University\\
         College Station\\
         TX \ 77843\\
         USA}
\email{sottile@math.tamu.edu}
\urladdr{www.math.tamu.edu/\~{}sottile}
\thanks{Work of Sottile supported by NSF grant DMS-0701050}
\subjclass[2000]{14M15, 14N15}
\keywords{Schubert calculus, isotropic Schubert calculus, transversality}
\begin{abstract}
 We show that general isotropic flags for odd-orthogonal and symplectic
 groups are general for Schubert calculus on the classical Grassmannian in that
 Schubert varieties defined by such flags meet transversally.
 This strengthens a result of Belkale and Kumar.
\end{abstract}
\maketitle
%
%
Schubert varieties \Blue{$\Omega_IE_\bullet$} in a classical flag manifold $G/P$ are given by
a flag $E_\bullet$ and a Schubert condition $I$~\cite{FuPr}.
By Kleiman's Transversality Theorem~\cite{Kl74}, if the flags
$E_\bullet^1,\dotsc,E_\bullet^s$ are general, then any corresponding Schubert
varieties intersect (generically) transversally in that they meet transversally along a
Zariski dense open subset of every component of their intersection.

Oftentimes we do not have the luxury of general flags, yet need to show that the 
Schubert varieties meet transversally.
It is often sufficient for their intersection to be \Blue{{\it proper}} (has the expected
dimension or is empty).   
Belkale and Kumar~\cite{BK} needed such a case where  
$G$ was either $Sp(2n)$ or $SO(2n{+}1)$, the flags $E_\bullet$ were isotropic flags,
and $G/P$ was an isotropic Grassmannian which is naturally a subset of a classical
Grassmannian \Blue{$Gr$}. 

\begin{proposition}[Belkale and Kumar~\cite{BK}]\label{Prop:BK}
 The intersection $\cap_{i=1}^s\Omega_{I^s}E_\bullet^i$ in $Gr$ is proper
 when $E_\bullet^1,\dotsc,E_\bullet^s$ are general isotropic flags and  $I^1,\dotsc,I^s$
 are Schubert conditions for $Sp(2n)$ or $SO(2n{+}1)$
\end{proposition}

We show that the intersection is in fact transverse.

\begin{theorem}\label{Th:transversality}
 The intersection $\cap_{i=1}^s\Omega_{I^i}E_\bullet^i$ in $Gr$ is
 transverse when $E_\bullet^1,\dotsc,E_\bullet^s$ are general isotropic flags for $Sp(2n)$
 or $SO(2n{+}1)$. 
\end{theorem}

We use a case for $Gr$ where the flags are not
general, yet the corresponding Schubert varieties meet transversally.
Let $f_1(t),\dotsc,f_m(t)$ be a basis for the space of polynomials of
degree less than $m$.
These define a \Blue{{\it rational normal curve}} $\Blue{\gamma}\colon \C\to \C^m$
by 
\[
   \gamma\ \colon\ t\ \longmapsto\ (f_1(t),\dotsc,f_m(t))^T\,.
\]
(We use column vectors and $(\dotsb)^T$ denotes transpose.)
For each $t\in \C$ and $i=1,\dotsc,m$, the $i$-plane 
osculating $\gamma$ at $\gamma(t)$ is the linear span of
$\gamma(t),\gamma'(t),\dotsc,\gamma^{(i-1)}(t)$.
These osculating planes form the 
\Blue{{\it osculating flag}} $\Blue{E_\bullet(t)}$.
An intersection of Schubert varieties for $Gr$ given by osculating flags
consists of linear series on ${\mathbb P}^1$ with at least some prescribed ramification.
Eisenbud and Harris~\cite{EH} showed that this intersection is proper.

\begin{proposition}\label{P:EH}
 The intersection $\cap_{i=1}^s\Omega_{I^i}E_\bullet(t_i)$ in $Gr$ is
 proper if $t_1,\dotsc,t_s\in\C$ are distinct.
\end{proposition}

This result is elementary---the codimension of the Schubert variety $\Omega_IE_\bullet(t)$
is the order of vanishing at $t$ of the Wronskian of a general linear series in
$\Omega_IE_\bullet(t)$.
Considerably less elementary is the following result of Mukhin, Tarasov, and
Varchenko~\cite[Corollary 6.3]{MTV}. 
 
\begin{proposition}\label{P:MTV}
 The intersection $\cap_{i=1}^s\Omega_{I^i}E_\bullet(t_i)$ in $Gr$ is
 transverse if $t_1,\dotsc,t_s\in\R$ are distinct.
\end{proposition}

Mukhin, Tarasov, and Varchenko proved this when the intersection is 
zero-dimensional, but the full statement follows
from their result via a standard argument.
Suppose that an intersection of Schubert varieties as in Proposition~\ref{P:MTV} has
dimension $r(>0)$ and let $Z$ be any of its components.
Let $\iota$ be the codimension 1 Schubert condition, so that $\Omega_\iota E_\bullet$ is a
hypersurface in $Gr$.
Let $u_1,\dotsc,u_r\in\R$ be distinct from $t_1,\dotsc,t_s$.
Then the intersection
 \begin{equation}\label{Eq:zero}
   \bigcap_{i=1}^s\Omega_{I^i}E_\bullet(t_i) \ \cap\ \bigcap_{i=1}^r\Omega_\iota E_\bullet(u_i)
 \end{equation}
is zero-dimensional and therefore transverse.
Since $\Omega_\iota E_\bullet$ meets every curve in $Gr$, the intersection
$Z\cap \cap_{i=1}^r\Omega_\iota E_\bullet(u_i)$ is non-empty.
Thus the intersection of Proposition~\ref{P:MTV} was transverse along $Z$, for
otherwise the intersection~\eqref{Eq:zero} would not be transverse at points of $Z$.

Let $\langle\,,\,\rangle$ be a non-degenerate alternating form on 
$\C^{2n}$ whose matrix $(\langle e_i,e_j\rangle)_{i,j=1,\dotsc,n}$
with respect to the standard ordered basis $e_1,\dotsc,e_{2n}$ is 
\[
    \left(\begin{array}{rr}0&J\\-J&0\end{array}\right)\ ,
\]
where $J$ is the anti-diagonal matrix $(1,\dotsc,1)$ of size $n$.
The \Blue{{\it symplectic group} $Sp(2n)$} is the group of linear transformations of
$\C^{2n}$ which preserve $\langle\,,\,\rangle$.
In this ordered basis 
 \begin{equation}\label{Eq:SP2n}
   \Blue{\gamma(t)}\ :=\ \left(
      1,\,t,\,\frac{t^2}{2},\,\ldots,\, \frac{t^n}{n!},\,
  -\frac{t^{n+1}}{(n+1)!},\,\frac{t^{n+2}}{(n+2)!},\,
      \ldots,\,(-1)^{n-1}\frac{t^{2n-1}}{(2n-1)!}\right)^T
 \end{equation}
is a rational normal curve whose osculating flag 
is \Blue{{\it isotropic}} in that $E_{2n-i}(t)$ annihilates $E_i(t)$ for $i< 2n$.
We leave this as an exercise for the reader.

Similarly, let $\langle\,,\,\rangle$  be a non-degenerate symmetric form 
on $\C^{2n+1}$ whose matrix is the
anti-diagonal matrix $(1,\dotsc,1)$ of size $2n+1$.
The \Blue{{\it special orthogonal group} $SO(2n{+}1)$} is the group of linear
transformations of $\C^{2n+1}$ of determinant 1 which preserve $\langle\,,\,\rangle$.
Then
 \begin{equation}\label{Eq:SO2n+1}
   \Blue{\gamma(t)}\ :=\ \left( 
      1,\,t,\,\frac{t^2}{2},\,\ldots,\, \frac{t^n}{n!},\,
     -\frac{t^{n+1}}{(n+1)!},\,\frac{t^{n+2}}{(n+2)!},\,
      \ldots,\,(-1)^n\frac{t^{2n}}{(2n)!}\right)^T
 \end{equation}
is a rational normal curve whose osculating flag 
is \Blue{{\it isotropic}} in that $E_{2n+1-i}(t)$ annihilates $E_i(t)$ for $i\leq 2n$.

Since it is an open condition on $s$-tuples of isotropic flags that Schubert varieties in $Gr$
meet properly or meet transversally, 
Proposition~\ref{Prop:BK} and Theorem~\ref{Th:transversality}
follow from Propositions~\ref{P:EH} and~\ref{P:MTV}, respectively.
These rational normal curves~\eqref{Eq:SP2n} and~\eqref{Eq:SO2n+1} were
introduced in~\cite{S} to study the analog of the Shapiro conjecture~\cite{So:Shapiro} 
for flag varieties for $Sp(2n)$ and $SO(2n{+}1)$, and the proof of the
Shapiro conjecture~\cite{MTV_Sh} motivated Proposition~\ref{P:MTV}.

These special osculating flags are better understood in terms of Lie theory.
Let $G$ be a semisimple complex Lie group with Lie algebra $\frakg$.
The adjoint action of $G$ on the nilpotent elements of $\frakg$ has finitely many orbits,
with dense orbit consisting of \Blue{{\it principal nilpotent}} elements of $\frakg$.
Write $\Blue{\exp}\colon \frakg\to G$ for the exponential map.
For a principal nilpotent $\eta\in\frakg$, 
$\{\exp(t\eta)\mid t\in \C\}$ is the corresponding 1-parameter
subgroup of $G$.  
It is natural to consider
Schubert varieties defined by translates
of a fixed flag by elements $\exp(t\eta)$.

The matrix $\eta\in\mathfrak{sl}_m$ with entries $1,2,\dotsc,m{-}1$ below its diagonal
is principal nilpotent.
Dale Peterson observed that 
the action  of $\exp(t\eta)$ on the standard coordinate flag  gives the osculating flag
$E_\bullet(t)$ to the rational normal curve $\gamma(t):=(1,t,t^2,\dotsc,t^{m-1})^T$.
The osculating flags to~\eqref{Eq:SP2n} and~\eqref{Eq:SO2n+1} also 
arise from exponentiating principal nilpotents in $\mathfrak{sp}_{2n}$ and 
$\mathfrak{so}_{2n+1}$, respectively.
These nilpotents have entries $1,\dotsc,1,-1,\dotsc,-1$ below their diagonals with $n$ 1s.
We obtain flags osculating a rational normal curve because principal nilpotents are mapped
to principal nilpotents under the inclusions
$\mathfrak{sp}_{2n}\hookrightarrow\mathfrak{sl}_{2n}$ and
$\mathfrak{so}_{2n+1}\hookrightarrow\mathfrak{sl}_{2n+1}$.

This is not the case for the even orthogonal groups, which explains their exclusion from
Theorem~\ref{Th:transversality}. 
A principal nilpotent for $\mathfrak{so}_{2n}$ is the $2n \times 2n$ matrix $\eta$ with
$1$ in positions $i,i{+}1$ and $-1$ in positions $2n{-}i,2n{-}i{+}1$ for $i=1,\dotsc,n$
(it has $1,\dotsc,1,0,-1,\dotsc,-1$ below its diagonal) and also 
$1$ in position $n{-}1,n{+}1$ and $-1$ in position $n,n{+}2$.
As  $\eta^{2n-1}=0$, it is not a principal nilpotent for $\mathfrak{sl}_{2n}$,
whose principal nilpotents $N$ have $N^{2n-1}\neq 0$.

We point out a further limitation of this method.
Proposition~\ref{P:EH} becomes false if we replace a classical Grassmannian $Gr$ by a
general type $A$ flag variety.
Indeed, in the 8-dimensional manifold of flags $\{F_1\subset F_3\subset\C^5\}$ consisting
of a 1-dimensional subspace lying in a 3-dimensional subspace in $\C^5$, the Schubert
variety $\Omega_{32514}E_\bullet$ has codimension 5 and the Schubert variety
$\Omega_{21435}E_\bullet$ has codimension 2.
Consequently, if $E_\bullet, E'_\bullet$, and $E''_\bullet$ are general flags, then
\[
   \Omega_{32514}E_\bullet\,\cap\, \Omega_{21435}E'_\bullet\,\cap\, \Omega_{21435}E''_\bullet
\]
is empty for dimension reasons.
If however, $E_\bullet, E'_\bullet$, and $E''_\bullet$ osculate a rational normal
curve, then the intersection is non-empty.
This is shown in Section 3.3.6 of~\cite{RSSS}.

\providecommand{\bysame}{\leavevmode\hbox to3em{\hrulefill}\thinspace}
\providecommand{\MR}{\relax\ifhmode\unskip\space\fi MR }
\providecommand{\MRhref}[2]{%
  \href{http://www.ams.org/mathscinet-getitem?mr=#1}{#2}
}
\providecommand{\href}[2]{#2}

\end{document}